\documentstyle[amssymb,amsfonts,psfig]{amsart}

\def\R{{\hbox{\bf R}}}
\def\C{{\hbox{\bf C}}}

 at 10 true pt

\def\ts{{T_\sigma}}

\def\allt#1{%
\smash{
 \vtop{%
     \ialign{%
        ##\crcr
        $\hfil\displaystyle{\tilde \forall}\hfil$\crcr%
        \noalign{\kern1.5pt\nointerlineskip}
        $\hfil\!\!#1\hfil$\crcr\noalign{\kern1.5pt}
        }
       }
      } \hbox{$\vphantom{#1}$}
     }

\def\be#1{\begin{equation}\label{#1}}
\def\bas{\begin{align*}}
\def\eas{\end{align*}}
\def\bi{\begin{itemize}}
\def\ei{\end{itemize}}

\def\eps{\varepsilon}

\def \endprf{\hfill  {\vrule height6pt width6pt depth0pt}\medskip}
\def\emph#1{{\it #1}}
\def\textbf#1{{\bf #1}}

\def\calA{{\cal A}}
\def\F{{\bf F}}

\theoremstyle{plain}
  
  \newtheorem{theo}[subsection]{Theorem}

\theoremstyle{remark}

\theoremstyle{definition}

\include{psfig}

\begin{document}

\title{On a problem by Arens, Goldberg, and Luxemburg}

\begin{abstract}
We construct a normed algebra ${\calA}$ 
with norm $N(\cdot)$ over the reals, which is {\em quadrative} in the sense 
that $N(A^2) \le N(A)^2$ for all $A \in {\cal A}$, but is not 3-{\em 
bounded} in the sense that $N(A^3) \le N(A)^3$.  This answers a question of 
Arens, Goldberg, and Luxemburg.
\end{abstract}

\author{Raymond Redheffer}
\address{Department of Mathematics, UCLA, Los Angeles CA 90095-1555}
\email{rr@@math.ucla.edu}

\author{Terence Tao}
\address{Department of Mathematics, UCLA, Los Angeles CA 90095-1555}
\email{tao@@math.ucla.edu}

\maketitle

Let $\calA$ be a normed algebra over a field $\F$, either $\R$ or
$\C$. In \cite{AGL} the norm $N$ of the algebra is called {\em quadrative}
if 
$$
N(A^2) \le N(A)^2 \quad {\rm for \; all} \quad A \in \calA,
$$
$k$-{\em bounded} for a positive integer $k$ if $N(A^k) \le N(A)^k$ for all
$A \in \calA$, and {\em strongly stable} if it is $k$-bounded for all $k =
1,2,3,\ldots$. It was seen in \cite{AG} that boundedness for a particular
$k > 2$ does not ensure strong stability or even quadrativity. Let $W =
(\omega_{ij})$ be a fixed $2\times 2$ matrix of positive entries. Then for the
$W$-{\em weighted} sup norm on $\C^{2\times 2}$, the algebra of $n\times n$
complex matrices,
$$
||A||_{W,\infty} = \max_{i,j} \, \omega_{ij} \, |\alpha_{ij}|, \qquad
A = (\alpha_{ij}) \in \C^{2\times 2},
$$
Arens and Goldberg proved:
\medskip

\noindent {\bf Theorem} [1, Theorem 2]. {\em If $k \ge 3$, then there
exists a $2\times 2$ weight matrix $W$ for which $||\cdot||_{W,\infty}$ is
$k$-bounded but not strongly stable, in fact not even quadrative on
$\C^{2\times 2}$.}
\medskip

Our main theorem gives a negative answer to the following question raised in 
\cite{AGL}: {\em Does quadrativity imply strong stability\/?}

\begin{theo} \label{th1}
There exists a commutative algebra $\calA$ of $2 \times 2$ matrices over $\R$
and a norm on $\calA$ such that $|A^2| \le |A|^2$ for all $A \in \calA$ and
$|A^3| > |A|^3$ for some $A \in \calA$.
\end{theo}

\noindent {\em Proof.} In Theorem \ref{th1} the elements of the algebra are
real matrices of the form
$$
A = \left(\begin{array}{ll}
a & b \\
0 & a
\end{array}\right) =: [[a,\,b]] 
$$
where the symbol on the right is introduced to save space. We will use the
identity
$$
[[a,b]]^k = [[a^k, k ba^{k-1}]]
$$
for any $[[a,b]]$ and any integer $k \ge 1$.

We observe that the algebra ${\calA}$ contains a multiplicative semi-group
$$
G := \{ \exp[[-t,t]]: t \ge 0 \} = \{ [[e^{-t}, t e^{-t}]]: t \ge 0\}
$$
In particular, if $A = [[e^{-t}, t e^{-t}]]$ is an element of $G$, then so is
$A^2 = [[e^{-2t}, 2t e^{-2t}]]$ and $A^3 = [[e^{-3t}, 3t e^{-3t}]]$.  

We can write $G$ as a graph of $b$ over $a$; indeed, setting $a = e^{-t}$ we 
have
$$
G := \{ [[a,b]]: 0 < a \le 1, b = f(a) \}
$$
where $f(a)$ is the function $f(a) := -a \log a$ on the interval $\{ 0 < a
\le 1 \}$. We remark that on this interval the function $f$ is concave (since
$f''(a) = -1/a$), non-negative and attains its maximum at the point $a = 
e^{-1}$, $f(a) = e^{-1}$. We define the modified function $g(a)$ on $\{ 0 < 
a \le 1 \}$ by setting $g(a) := f(a)$ when $e^{-1} \le a \le 1$ and $g(a) := 
e^{-1}$ when $0 < a \le e^{-1}$; note that $g$ is still (weakly) concave.

Define a {\em ball} to be any non-empty bounded open convex subset of $\calA$
which is symmetric around the origin. Then for every ball $\Omega$, we can
define a norm $N_\Omega$ on $\calA$ in the usual manner as
$$
N_\Omega(A) := \inf \{ t: t > 0, A \in t \Omega \},
$$
so that $\Omega$ is the unit ball of $A$. The fact that $\Omega$ is a ball
ensures that $N_\Omega$ is indeed a norm.

Let $k \ge 2$ be an integer. We say that a norm $N(\cdot)$ on $\calA$ is
$k$-{\em bounded} if one has $N(A^k) \le N(A)^k$ for all $A \in \calA$. Also,
we shall say that a ball $\Omega$ is $k$-{\em bounded} if one has $A^k \in
\Omega$ whenever $A \in \Omega$. It is clear from homogeneity that $N_\Omega$
is $k$-bounded if and only if $\Omega$ is $k$-bounded. We say that $N$ or
$\Omega$ is {\em quadrative} if it is $2$-bounded.  

As an example, consider the set
$$
\Omega_0 := \{ [[a,b]]: |a| < 1; |b| < g(|a|) \}.
$$
It is clear that this set is a ball. We now show that $\Omega_0$ is $k$-bounded
for every integer $k \ge 2$. Let $[[a,b]] \in \Omega_0$; we have to show that
$[[a,b]]^k = [[a^k, k ba^{k-1}]]$ is also in $\Omega_0$.

\newpage

\begin{figure}[htbp] 
\centering
\ \psfig{figure=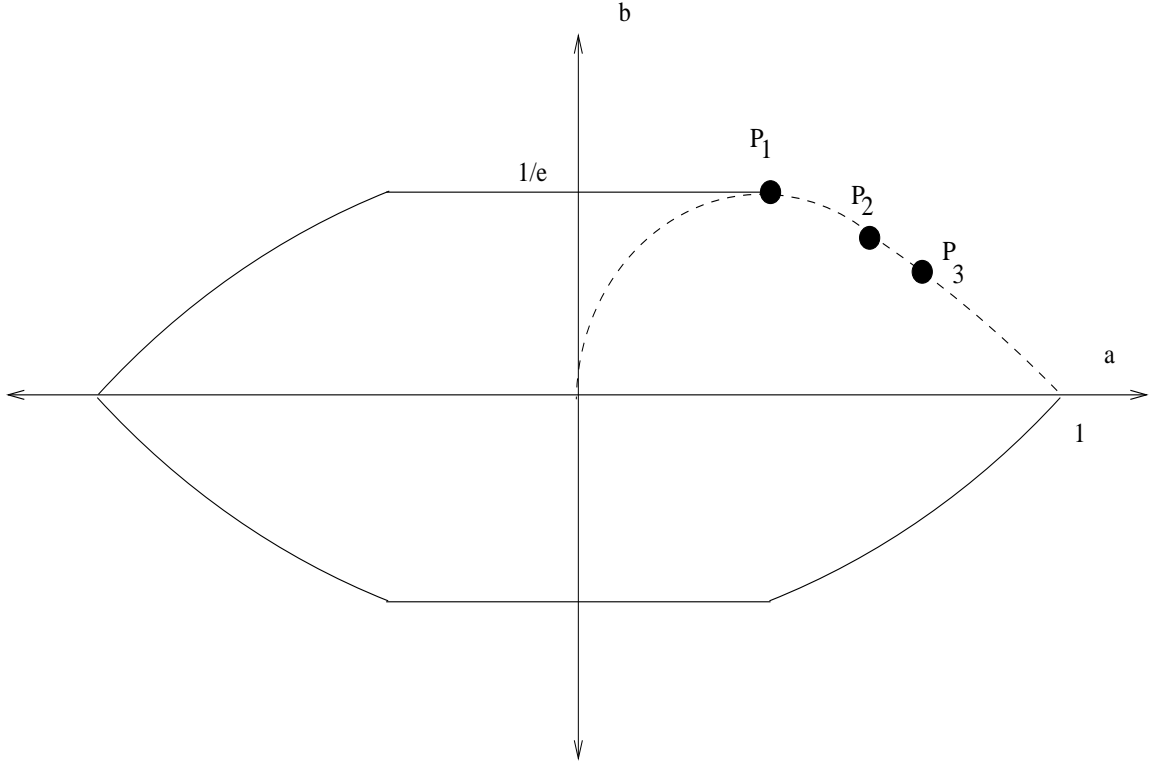,height=4in,width=6in} 
\caption{The ball $\Omega_0$, together with the three points $P_1$, $P_2$, $P_3$.  The dotted curve is the graph of $b = f(a)$
for $0 < a \le 1$; this is the semi-group $G$.}
\end{figure}

By reflection symmetry in the $a$ and $b$ axes, we may assume that we are in
the first quadrant $a, b \ge 0$. There are two cases: $e^{-1} \le a < 1$ and 
$0
\le a \le e^{-1}$.

First suppose that $e^{-1} \le a < 1$. Then $b < g(a) = -a \log a$. Thus
$$
k b a^{k-1} < -a^k \log a^k = f(a^k) \le g(a^k),
$$
and so $[[a^k, k ba^{k-1}]] \in \Omega_0$ as desired.

Now suppose that $0 \le a \le e^{-1}$. Then $b < g(a) = e^{-1}$.  Thus
$$
k b a^{k-1} < k e^{-k} \le e^{-1} = g(a^k)
$$
since the function $t e^{-t}$ attains its maximum at $t=1$, and since $a^k$ is
clearly bounded by $e^{-1}$. Thus $[[a^k, k ba^{k-1}]] \in \Omega_0$ as 
desired.

We identify three interesting points on the boundary of $\Omega_0$: $P_1 :=
[[e^{-1}, e^{-1}]]$, $P_2 := [[e^{-1/2}, \frac{1}{2} e^{-1/2}]]$, and $P_3 :=
[[e^{-1/3}, \frac{1}{3} e^{-1/3}]]$. Note that $P_3^3 = P_1$ and $P_2^2 = P_1$.
Also, $P_1$ is the point of $\Omega_0$ where the two constraints $|b| < 
f(|a|)$ and $|b| < e^{-1}$ intersect.

We now modify the ball $\Omega_0$ slightly, to prove
\medskip

\noindent {\bf Proposition.} {\em There exists a ball $\Omega$ which is 
$2$-bounded but not $3$-bounded.}
\medskip

\noindent {\em Proof.} The idea is to chip a small amount away from
$\Omega_0$, enough to destroy the 3-boundedness but not enough to destroy
the 2-boundedness.

We shall need three small numbers $0 < \eps_3 \ll \eps_2 \ll \eps_1 < 1$ to
be chosen later. We define $\Omega$ to be the set of matrices $[[a,b]]$ in
which $|a| < e^{-1/3}$ and $b$ satisfies all three of the inequalities
$$
|b| < g(|a|), \quad |b| < e^{-1} - \eps_3, \quad |b| < e^{-1/2} - 
{\ts \frac{1}{2}} |a| - \eps_1.
$$
Note that the line $b = e^{-1/2} - \frac{1}{2} a$ is the tangent line to the
curve $b = g(a)$ at the point $P_2$. Thus the restriction $|b| \le e^{-1/2}
- \frac{1}{2} |a| - \eps_1$ cuts off a small sliver of $\Omega_0$ near the
point $P_2$ (and similarly for the other three quadrants, by reflection
symmetry). The restriction $|a| < e^{-1/3}$ cuts off everything in 
$\Omega_0$ to the right of $P_3$, while the restriction $|b| < e^{-1} -
\eps_3$ cuts off a very thin horizontal sliver from the straight portion of 
the boundary of $\Omega_0$, and in particular cuts off a small sliver near 
$P_1$.

\begin{figure}[htbp]
\centering
\ \psfig{figure=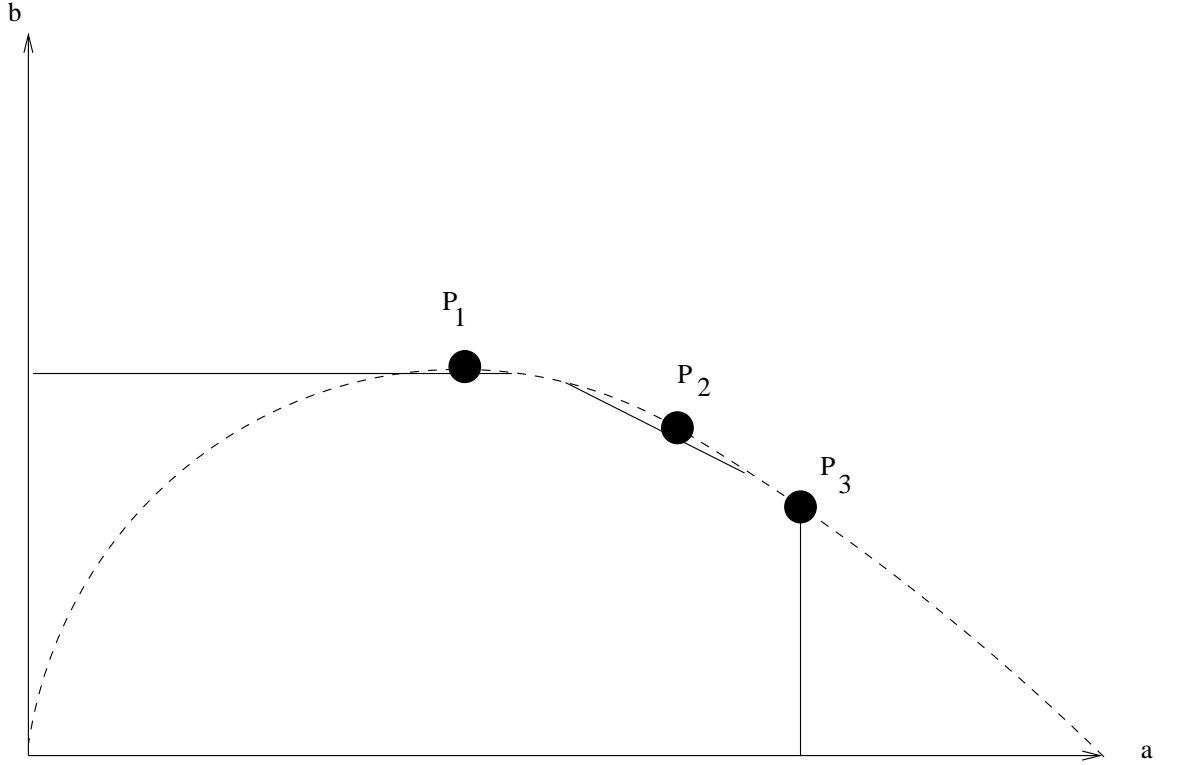,height=4in,width=6in}
\caption{The ball $\Omega$ in the first quadrant. The vertical line is the 
condition $|a| < e^{-1/3}$; the horizontal line is the condition $|b| < 
e^{-1} - \eps_3$; and the slanted line is the condition $|b| < e^{-1/2} 
- \frac{1}{2} |a| - \eps_1$.  Note that $\Omega$ gets arbitrarily close to 
$P_3$ but not to $P_1$ or $P_2$; also the region removed near $P_2$ is 
larger than that near $P_1$ since it depends on $\eps_1$ instead of 
$\eps_3$.}
\end{figure}

It is clear that $\Omega$ is still a bounded open non-empty convex symmetric 
set, i.e. a ball.  Also, it is clear that $\Omega$ is no longer 3-bounded, 
because one can get arbitrarily close to $P_3$ in $\Omega$, but one cannot
get arbitrarily close to $P_1 = P_3^3$.

It remains to show that $\Omega$ is 2-bounded. To do this, we take any
$[[a,b]] \in \Omega$; our task is to show that $[[a,b]]^2 = [[a^2, 2ab]]$ is 
also in $\Omega$. By symmetry we may assume that we are in the first 
quadrant $a,b \ge 0$. Since $a < e^{-1/3}$, we have $a^2 < e^{-2/3} < 
e^{-1/3}$, so we only have to show the three inequalities
\begin{eqnarray}
2ab & < & g(a^2)\label{a-1}\\
2ab & < & e^{-1} - \eps_3 \label{a-2}\\
2ab & < & e^{-1/2} - {\ts\frac{1}{2}} a^2 - \eps_1.\label{a-3}
\end{eqnarray}
Recall that the line $y = e^{-1/2} - \frac{1}{2} x - \eps_1$ was just a small
perturbation of the tangent line $y = e^{-1/2} - \frac{1}{2} x$ of the 
curve $y=g(x)$ at the point $x = e^{-1/2}, y = \frac{1}{2} e^{-1/2}$. In 
particular we see from the concavity of $g$ that, if $\eps_1$ is 
sufficiently small,
$$
e^{-1/2} - {\ts\frac{1}{2}} x - \eps_1 > g(x)
$$
for all $x < e^{-2/3} < e^{-1/2}$.  Since $a^2 < e^{-2/3}$, we thus see that
the condition (\ref{a-3}) is redundant, being implied automatically by 
(\ref{a-1}).

It remains to prove (\ref{a-1}) and (\ref{a-2}). To do this we divide into
three cases.
\medskip

\noindent {\bf Case 1.} $e^{-1/2}+\eps_2 \le a < e^{-1/3}$. Then we have
$b < g(a) = f(a) = -a \log a$.  Thus $$ 2ab < -a^2 \log a^2 = f(a^2) = 
g(a^2)$$
since
\begin{equation}\label{a-bound}
e^{-1} + 2 \eps_2 e^{-1/2} + \eps_2^2 \le a^2 < e^{-2/3}.
\end{equation}
This gives (\ref{a-1}). If $\eps_2$ is chosen sufficiently small compared
to $\eps_1$, and $\eps_3$ is chosen sufficiently small compared to $\eps_2$, 
then we see from (\ref{a-bound}) (and the monotonicity of $g(x)$ for $x > 
e^{-1}$) that
$$
g(a^2) < g(e^{-1}) - \eps_3 = e^{-1} - \eps_3.
$$
This gives (\ref{a-2}) as desired.
\medskip

\noindent {\bf Case 2.} $e^{-1/2} - \eps_2 < a < e^{-1/2} + \eps_2$. In
this case we use the bound
$$
b < e^{-1/2} - {\ts\frac{1}{2}} a - \eps_1,
$$
and so
$$
2ab < 2 e^{-1/2} a - a^2 - 2 \eps_1 a.
$$
Since $a = e^{-1/2} + O(\eps_2)$, we thus have
$$
2ab < e^{-1} - 2 \eps_1 e^{-1/2} + O(\eps_2).
$$
On the other hand, since $a^2 = e^{-1} + O(\eps_2)$, we have
$$
g(a^2) = g(e^{-1}) + O(\eps_2) = e^{-1} + O(\eps_2).
$$
Thus if $\eps_2$ is sufficiently small compared to $\eps_1$, we obtain 
(\ref{a-1}). Using the above estimate for $2ab$, we see that (\ref{a-2})
follows from
$$
e^{-1} - 2\eps_1 e^{-1/2} + O(\eps_2) < e^{-1} - \eps_3.
$$
This holds if both $\eps_2$ and $\eps_3$ are sufficiently small compared to
$\eps_1$.
\medskip

\noindent {\bf Case 3.} $0 < a \le e^{-1/2} - \eps_2$. In this case we use 
the bound $b < g(a)$, so that $2ab < 2a g(a)$. Since
\begin{equation}\label{ab-2}
a^2 \le e^{-1} - 2 \eps_2 e^{-1/2} + \eps_2^2
\end{equation}
we have $g(a^2) = e^{-1}$ where $\eps_2$ is small. Thus (\ref{a-1}) follows 
from (\ref{a-2}), and it suffices to show that
$$
2 a g(a) < e^{-1} - \eps_3.
$$
First suppose that $a \le e^{-1}$. Then $2 a g(a) \le 2 e^{-2}$, which is 
certainly acceptable if $\eps_3$ is small enough. Thus we may take $a > 
e^{-1}$, in which case
$$
2 a g(a) = 2 a f(a) = - a^2 \log a^2 = f(a^2).
$$
Since $f$ attains its maximum $e^{-1}$ at $e^{-1}$, we thus see from 
(\ref{ab-2}) that $f(a^2) < e^{-1} - \eps_3$, if $\eps_3$ is sufficiently 
small compared to $\eps_2$. This concludes the proof of the Proposition 
in all three cases. \hfill $\square$
\bigskip

One may try to improve this counterexample by adding another natural 
condition to the norm $N$, namely that the identity $[[1,0]]$ have norm 1.  
This is equivalent to $[[1,0]]$ lying on the boundary of $\Omega$. It is 
true that the counterexample constructed above does not obey this condition, 
but this is easily rectified by replacing the ball $\Omega$ constructed 
above with the convex hull ${\rm hull}(\Omega, [[1,0]], [[-1,0]])$ of 
$\Omega$ with the points $[[\pm 1, 0]]$.  We omit the computation which 
shows that this ball remains 2-bounded and not 3-bounded.
\bigskip

\noindent {\bf Acknowledgement.} It is a a pleasure to express our 
appreciation
to Professor Moshe Goldberg, who brought the problem to our attention and supplied
the references.

\end{document}